\documentclass[12pt]{article}
\usepackage{amssymb,amsmath}
\usepackage{hyperref}
\usepackage{graphicx}
\usepackage{color}
\usepackage{chngcntr}

\begin{document}

%\begin{center}
\title{\LARGE\bf The Fourier expansion approximation for high-accuracy computation of the Voigt/complex error function at \\
small imaginary argument}

%\bigskip
\author{\normalsize\bf S. M. Abrarov\footnote{\scriptsize{Dept. Earth and Space Science and Engineering, York University, Toronto, Canada, M3J 1P3.}}\, and B. M. Quine$^{*}$\footnote{\scriptsize{Dept. Physics and Astronomy, York University, Toronto, Canada, M3J 1P3.}}}

\date{June 25, 2016}
\maketitle
%\vspace{1cm}%\bigskip

\begin{abstract}
It is known that the computation of the Voigt/complex error function is problematic for highly accurate and rapid computation at small imaginary argument $y << 1$, where $y = \operatorname{Im} \left[ z \right]$. In this paper we consider an approximation based on the Fourier expansion that can be used to resolve effectively such a problem when $y \to 0$.\\

\noindent {\bf Keywords:} Complex error function; Voigt function; Faddeeva function; complex probability function; spectral line broadening
\end{abstract}

The complex error function also known as the Faddeeva function is given by \cite{Hui1978, Weideman1994, Schreier1992, Zaghloul2011, Abrarov2011}
\begin{equation}\label{eq_1}
w\left( z \right)={{e}^{-{{z}^{2}}}}\left( 1+\frac{2i}{\sqrt{\pi }}\int\limits_{0}^{z}{{{e}^{{{t}^{2}}}}dt} \right),
\end{equation}
where $z = x + iy$ is the complex argument. The real part of the complex error function is known as the Voigt function \cite{Hui1978, Schreier1992, Armstrong1967, Letchworth2007} that are widely used in many fields of Physics, Chemistry and Astronomy. The integral on the right side of equation \eqref{eq_1} cannot be taken analytically. Consequently, the complex error function requires a numerical solution.

We have shown previously that using the Fourier expansion method the complex error function can be represented as \cite{Abrarov2011, Abrarov2012}
\footnotesize
\begin{equation}\label{eq_2}
\begin{aligned}
w\left( z \right) \approx & \frac{i}{{2\sqrt \pi  }}\left[ {\sum\limits_{n = 0}^N {{a_n}{\tau _m}\left( {\frac{{1 - {e^{i\,\,\left( {n\pi  + {\tau _m}z} \right)}}}}{{n\,\pi  + {\tau _m}z}} - \frac{{1 - {e^{i\,\,\left( { - n\pi  + {\tau _m}z} \right)}}}}{{n\,\pi  - {\tau _m}z}}} \right)}  - {a_0}\frac{{1 - {e^{i{\tau _m}z}}}}{z}} \right] \\
= & i\frac{{1 - {e^{i{\tau _m}z}}}}{{{\tau _m}z}} + i\frac{{\tau _m^2z}}{{\sqrt \pi  }}\sum\limits_{n = 1}^N {{a_n}\frac{{{{\left( { - 1} \right)}^n}{e^{i{\tau _m}z}} - 1}}{{{n^2}\,{\pi ^2} - \tau _m^2{z^2}}}}.
\end{aligned}
\end{equation}
\normalsize
where
$$
{a_n} \approx \frac{{2\sqrt \pi  }}{{{\tau _m}}}\exp \left( { - \frac{{{n^2}{\pi ^2}}}{{\tau _m^2}}} \right)
$$
is a set of the Fourier expansion coefficients, ${\tau_m}= 12$ and $N = 23$ (for simplicity we assume that $x = \operatorname{Re\left[z\right]}$ and $y = \operatorname{Im\left[z\right]}$ are non-negative).

The approximation \eqref{eq_2} is efficient in computation. In particular, this approximation has been implemented in the latest version of the \href{https://arxiv.org/abs/1407.0748}{{\it{RooFit}}} package for highly accurate and rapid computation of the Faddeeva function (C++ source code is provided in \cite{Karbach2014}). Furthermore, the equation \eqref{eq_2} alone can cover with high-accuracy the entire domain $0 < x < 40,000$ and $10^{-4} < y < 10^{2}$ required for the HITRAN molecular spectroscopic database \cite{Rothman2013}. However, as the input parameter $y$ decreases and becomes below  $10 ^{-6}$, the accuracy of the approximation \eqref{eq_2} deteriorates. Although the deterioration in accuracy of approximation \eqref{eq_2} with decreasing $y$ is much slower as compared to the Weideman\text{'}s approximation \cite{Weideman1994} (see also \cite{Abrarov2011} for more information), it cannot cover the narrow band domain along $x$-axis at $y << 1$. It should be noted that the computation of the Voigt/complex error function at small $y<<1$ is a common problem in many approximations \cite{Armstrong1967, Amamou2013, Abrarov2015}.

Let us show how a trivial rearrangement of the approximation \eqref{eq_2} can effectively resolve this problem. Consider the following identity \cite{Hui1978, Zaghloul2011}
$$
w\left(z \right)=2{{e}^{-{{z}^{2}}}}-w\left(-z \right)
$$
that can be simply rewritten as
\begin{equation}\label{eq_3}
w\left( z \right)={{e}^{-{{z}^{2}}}}+\frac{w\left( z \right)-w\left( -z \right)}{2}.
\end{equation}
Substituting the approximation \eqref{eq_2} into the right side of identity \eqref{eq_3} immediately leads to
\begin{equation}\label{eq_4}
w\left( z \right) \approx {{e}^{-{{z}^{2}}}}-i\frac{\cos \left( {{\tau }_{m}}z \right)-1}{{{\tau }_{m}}z}+i\frac{\tau _{m}^{2}z}{\sqrt{\pi }}\sum\limits_{n=1}^{N}{{{a}_{n}}\frac{{{\left( -1 \right)}^{n}}\cos \left( {{\tau }_{m}}z \right)-1}{{{n}^{2}}{{\pi }^{2}}-\tau _{m}^{2}{{z}^{2}}}}.
\end{equation}

Theoretically, the equations \eqref{eq_2} and \eqref{eq_4} are identical. Practically, however, these equations cover different domains. In particular, while the approximation \eqref{eq_2} provides high-accuracy at $y \gtrsim 10^{-6}$, the approximation \eqref{eq_4} is highly accurate at $y \lesssim 0.1$.

In order to quantify accuracy of the approximation \eqref{eq_4} we can define the relative errors for the real and imaginary parts as given by
$$
{{\Delta }_{\operatorname{Re}}}=\left| \frac{\operatorname{Re}\left[ {{w}_{\text{ref}\text{.}}}\left( z \right) \right]-\operatorname{Re}\left[ w\left( z \right) \right]}{\operatorname{Re}\left[ {{w}_{\text{ref}\text{.}}}\left( z \right) \right]} \right|
$$ 
and
$$
{{\Delta }_{\operatorname{Im}}}=\left| \frac{\operatorname{Im}\left[ {{w}_{\text{ref}\text{.}}}\left( z \right) \right]-\operatorname{Im}\left[ w\left( z \right) \right]}{\operatorname{Im}\left[ {{w}_{\text{ref}\text{.}}}\left( z \right) \right]} \right|,
$$
respectively, where $w_{{\rm{ref.}}} \left(z\right)$ is the reference. We obtained the highly accurate reference values by using Wolfram Mathematica (version 9) in enhanced precision mode.

%\newpage
\begin{figure}[ht]
\begin{center}
\includegraphics[width=24pc]{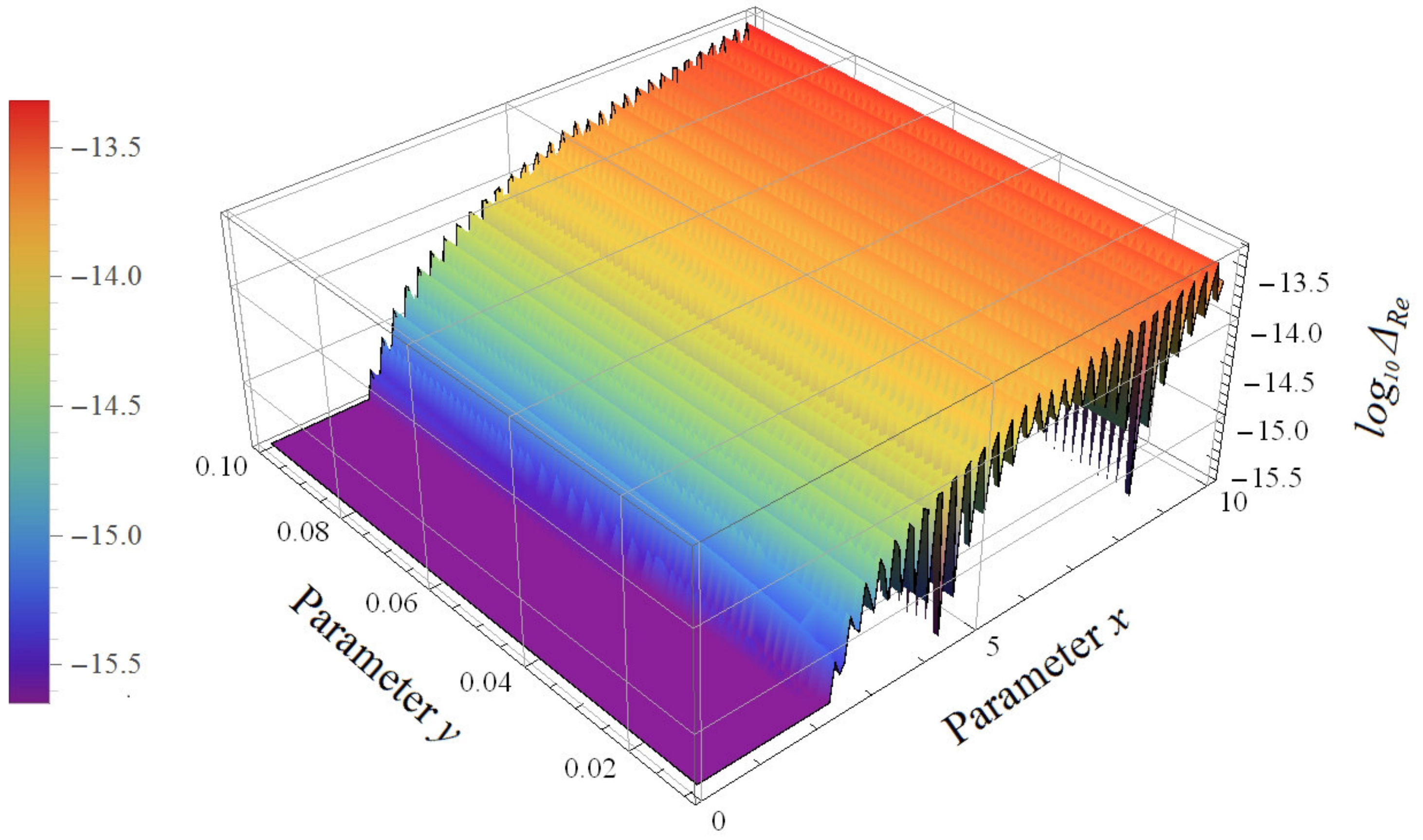}\hspace{1pc}%
\begin{minipage}[b]{26pc}
\vspace{0.15cm}
\small
{\sffamily {\bf{Fig. 1.}} Logarithm of the relative error $\log_{10}\Delta_{\operatorname{Re}}$ for the real part of the complex error function at $y \le 0.1$.}
\normalsize
\end{minipage}
\end{center}
\end{figure}

%\newpage
\begin{figure}[ht]
\begin{center}
\includegraphics[width=24pc]{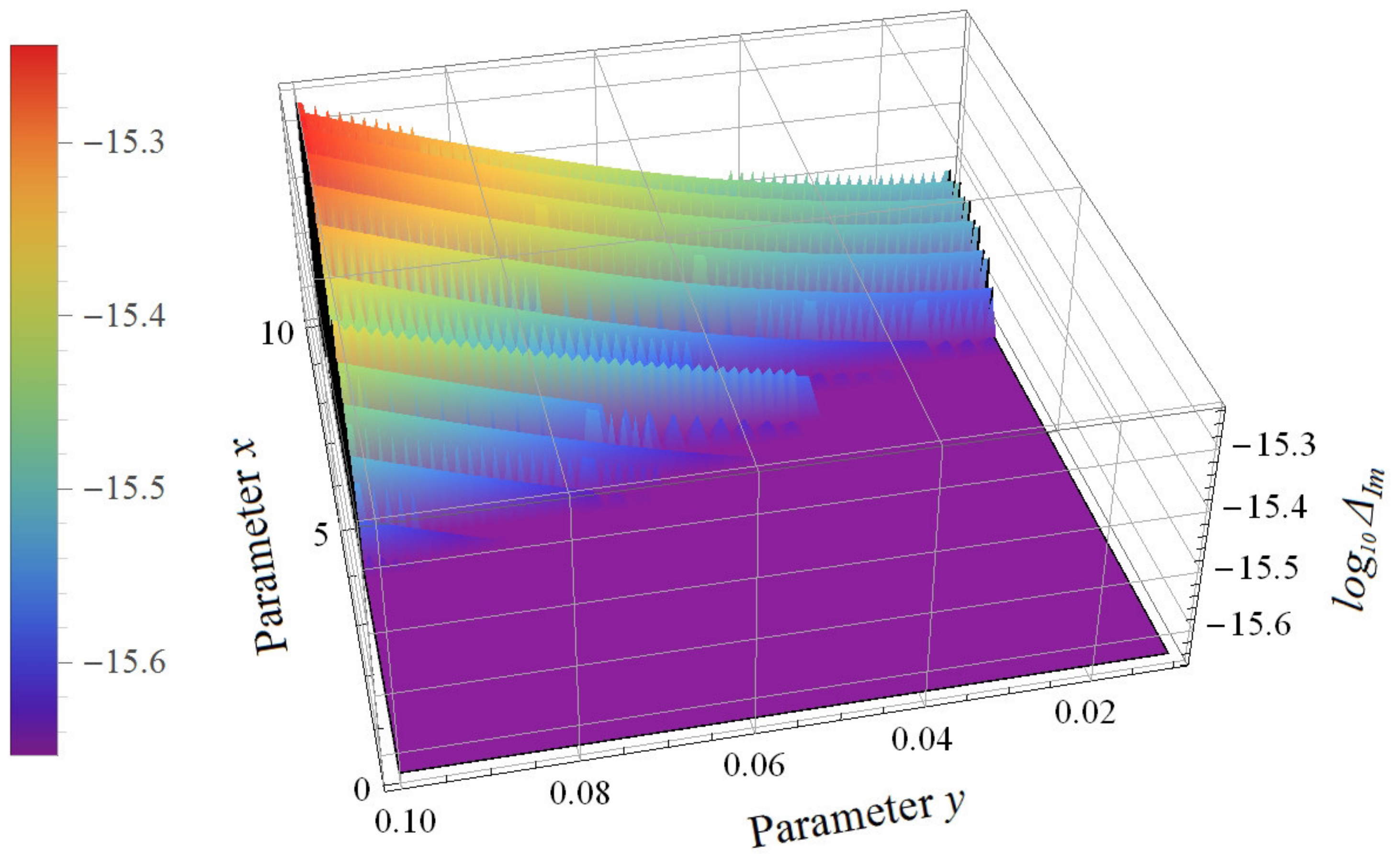}\hspace{1pc}%
\begin{minipage}[b]{26pc}
\vspace{0.15cm}
\small
{\sffamily {\bf{Fig. 2.}} Logarithm of the relative error $\log_{10}\Delta_{\operatorname{Im}}$ for the imaginary part of the complex error function at $y \le 0.1$.}
\normalsize
\end{minipage}
\end{center}
\end{figure}

Figure 1 shows the logarithm of the relative error for the real part of the complex error function. As we can see from this figure the worst accuracy is $\sim 10^{-13}$. The accuracy rapidly improves with decreasing $x$. Specifically, at $x \lesssim 2$ it becomes $\sim 10^{-16}$ as illustrated in Fig. 1 by dark-blue/purple area along $y$-axis.

Figure 2 depicts the logarithm of the relative error for the imaginary part of the complex error function. We can see from Fig. 2 that the accuracy is nearly perfect as the worst accuracy in the imaginary part is $\sim 10^{-15}$.

The computational test we performed shows that the approximation \eqref{eq_4} provides high-accuracy within the domain $0 \le x \le 10$ and $10^{-14} \le y \le 10^{-1}$. Thus, combining together the equations \eqref{eq_2} and \eqref{eq_4} we have developed an algorithm that covers with high-accuracy the entire domain required for practical applications.


\bigskip
\begin{thebibliography}{26}

\smallskip
\bibitem{Hui1978}
A.K. Hui, B.H. Armstrong and A.A. Wray, Rapid computation of the Voigt and complex error functions, J. Quantit. Spectrosc. Radiat. Transfer, 19 (5) (1978) 509-516. \\
\url{http://dx.doi.org/10.1016/0022-4073(78)90019-5}

\smallskip
\bibitem{Weideman1994}
J.A.C. Weideman, Computation of the complex error function, SIAM J. Numerical Analysis, 31 (1994) 1497-1518. \\
\url{http://dx.doi.org/10.1137/0731077}

\smallskip
\bibitem{Schreier1992}
F. Schreier, The Voigt and complex error function: A comparison of computational methods. J. Quant. Spectrosc. Radiat. Transfer, 48 (1992)
743-762. \\
\url{http://dx.doi.org/10.1016/0022-4073(92)90139-U}

\smallskip
\bibitem{Zaghloul2011}
M.R. Zaghloul and A.N. Ali, Algorithm 916: computing the Faddeyeva and Voigt functions, ACM Trans. Math. Software 38 (2011) 15:1-15:22. \\
\url{http://dx.doi.org/10.1145/2049673.2049679}

\smallskip
\bibitem{Abrarov2011}
S.M. Abrarov and B.M. Quine, Efficient algorithmic implementation of the Voigt/complex 	error function based on exponential series approximation, Appl. Math. Comp., 218 (2011) 1894-1902. \\
\url{http://dx.doi.org/10.1016/j.amc.2011.06.072}

\smallskip
\bibitem{Armstrong1967}
B.H. Armstrong, Spectrum line profilles: the Voigt function, J. Quant. Spectrosc. Radiat. Transfer. 7 (1967) 61-88. \\
\url{http://dx.doi.org/10.1016/0022-4073(67)90057-X}

\smallskip
\bibitem{Letchworth2007}
K.L. Letchworth and D.C. Benner, Rapid and accurate calculation of the Voigt function, J. Quant. Spectrosc. Radiat. Transfer, 107 (2007) 173-192. \\
\url{http://dx.doi.org/10.1016/j.jqsrt.2007.01.052}

\smallskip
\bibitem{Abrarov2012}
S.M. Abrarov and B.M. Quine, On the Fourier expansion method for highly accurate computation of the Voigt/complex error function in a rapid algorithm, \href{http://arxiv.org/abs/1205.1768}{arXiv:1205.1768} (2012).

\smallskip
\bibitem{Karbach2014}
T.M. Karbach, G. Raven and M. Schiller, Decay time integrals in neutral meson mixing and their efficient evaluation, \href{https://arxiv.org/abs/1407.0748}{arXiv:1407.0748} (2014).

\smallskip
\bibitem{Rothman2013}
L.S. Rothman, I.E. Gordon, Y. Babikov, A. Barbe, D.C. Benner, P.F. Bernath, M. Birk, L. Bizzocchi, V. Boudon, L.R. Brown, A. Campargue,
K. Chance, E.A. Cohen, L.H. Coudert, V.M. Devi, B.J. Drouin, A. Fayt, J.-M. Flaud, R.R. Gamache, J.J. Harrison, J.-M. Hartmann, C. Hill,
J.T. Hodges, D. Jacquemart, A. Jolly, J. Lamouroux, R.J. Le Roy, G. Li, D.A. Long, O.M. Lyulin, C.J. Mackie, S.T. Massie, S. Mikhailenko,
H.S.P. M\"uler, O.V. Naumenko, A.V. Nikitin, J. Orphal, V. Perevalov, A. Perrin, E.R. Polovtseva and C. Richard, The HITRAN2012 molecular
spectroscopic database, J. Quant. Spectrosc. Radiat. Transfer, 130 (2013) 4-50. \\
\url{http://dx.doi.org/10.1016/j.jqsrt.2013.07.002}

\smallskip
\bibitem{Amamou2013}
H. Amamou, B. Ferhat and A. Bois, Calculation of the Voigt function in the region of very small values of the parameter a where the calculation is notoriously difficult, Amer. J. Anal. Chem., 4 (2013) 725-731. \\
\url{http://dx.doi.org/10.4236/ajac.2013.412087}

\smallskip
\bibitem{Abrarov2015}
S.M. Abrarov and B.M. Quine, Accurate approximations for the complex error function with small imaginary argument, J. Math. Research, 7 (1) (2015) 44-53. \\
\url{http://dx.doi.org/10.5539/jmr.v1n1p44}

\end{thebibliography}
\end{document}